\documentclass[journal,10pt]{IEEEtran}
\usepackage{graphicx,epstopdf,multirow,multicol,booktabs,pgfplots}

\usepackage{amsmath,bm}
\usepackage{amssymb,amsthm}
\usepackage[T1]{fontenc}
\usepackage[colorlinks,citecolor=blue]{hyperref}
\allowdisplaybreaks
\usepackage{setspace}
\usepackage{cite}
\usepackage{siunitx}
\usepackage{times}
\usepackage{float}

\usepackage{amssymb}
\usepackage{amsmath}
\usepackage{accents}

\usepackage{xcolor}
\newcommand{\edit}[1]{\textcolor{black}{#1}}

\theoremstyle{remark}

\pgfplotsset{compat=1.14}

\begin{document}
\title{Operation of Natural Gas Pipeline Networks with Storage under Transient Flow Conditions}

\author{Sai Krishna Kanth Hari$^{\dagger}$, Kaarthik Sundar$^{*}$, Shriram Srinivasan$^{\dagger}$, Anatoly Zlotnik$^{\dagger}$, Russell Bent$^{\dagger}$
\thanks{$^{\dagger}$Applied Mathematics and Plasma Physics Group, Los Alamos National
Laboratory, Los Alamos, New Mexico, USA. E-mail: \texttt{\{hskkanth,shrirams,azlotnik,rbent\}@lanl.gov}}\;
\thanks{$^{*}$Information Systems and Modeling Group, Los Alamos National
Laboratory, Los Alamos, New Mexico, USA. E-mail: \texttt{{kaarthik}@lanl.gov}}\;
}

\maketitle

\begin{abstract}
We formulate a nonlinear optimal control problem for intra-day operation of a natural gas pipeline network that includes storage reservoirs. The dynamics of compressible gas flow through pipes, compressors, reservoirs, and wells are considered.  In particular, a reservoir is modeled as a rigid, hollow container that stores gas under isothermal conditions and uniform density, and a well is modeled as a vertical pipe. For each pipe, flow dynamics are described by a coupled partial differential equation (PDE) system in density and mass flux variables, with momentum dissipation modeled using the Darcy-Wiesbach friction approximation. Compressors are modeled as scaling up the pressure of gas between inlet and outlet. The governing equations for all network components are spatially discretized and assembled into a nonlinear differential-algebraic  equation (DAE) system, which synthesizes above-ground pipeline and subsurface reservoir dynamics into a single reduced-order model.  We seek to maximize an objective function that quantifies economic profit and network efficiency subject to the flow equations and inequalities that represent operating limitations.  The problem is solved using a primal-dual interior point solver and the solutions are validated in computational experiments and simulations on several  pipeline test networks to demonstrate the effectiveness of the proposed methodology. 
\end{abstract}

\IEEEpeerreviewmaketitle

\section{Introduction} \label{sec:intro}
Growth in natural gas usage worldwide is driven by increasing natural gas-fired electricity generation and lower natural gas prices.  In the United States, natural gas consumption reached a record high of 85 billion cubic feet per day in 2019 \cite{EIA-1}. Furthermore, data from the U.S. Energy Information Administration indicates that natural gas accounts for the largest share of the generation fuel mix since surpassing coal on an annual basis in 2016 \cite{EIA-2}. This increase presents several new challenges in gas pipeline system management.

First, new natural gas-fired plants are used to provide base load in power grid operations, and are also used to respond quickly to balance out fluctuations in electricity production by uncontrollable renewable sources \cite{Rinaldi2001,Li2008}. Pipeline loads are therefore subject to decision making for power grid operations, which has resulted in increasing and also more variable gas consumption, in contrast to historical daily consumption profiles that were steadier and more predictable. It is of growing interest to move beyond the traditional use of steady-state flow capacity models for managing daily gas pipeline markets and operations to more effectively utilize pipeline capacity given time-varying consumption of natural gas, by modeling the pipeline flows under transient conditions. Thus, developing models for natural gas storage facilities and incorporating them into the day-to-day operational problem is critical to enable pipeline management that is responsive to fluctuations caused by uncontrollable renewable sources of energy across all seasons and conditions. \edit{An interested reader is referred to \cite{Roald2020} for a framework that describes how to handle uncertainty in a joint gas-electric grid. }

The second challenge stems from the fact that weather is the primary driver of annual and seasonal fluctuations in natural gas consumption. During the winter season, U.S. natural gas consumption levels are at their highest because natural gas is the predominant fuel for residential and commercial heating. Accordingly, some pipeline systems are designed to include underground storage facilities, which are typically depleted petroleum reservoirs, that are necessary for these systems to function as intended.  Natural gas is injected into such facilities during the summer months, when supply is in excess of demand, and withdrawn from storage in the winter months when peak gas consumption exceeds the gas that can be transported into the system as flowing supply.  However, the nuanced relation between pressure gradients and flows in pipeline systems often preclude storage facilities and flowing supplies from being utilized at maximal capacities simultaneously.  Both of the above challenges can be addressed by joint optimal control of gas pipeline flows together with the supporting storage facilities.

The modeling of transient flow of natural gas through pipeline systems is a developed field, and was considered recently in the contexts of optimal control as well as state and parameter estimation \cite{Zlotnik2015,Sundar2018,Sundar2019,Jalving2018,Zlotnik2019}. In these studies, transient gas pipeline flow is modeled using simplifications of the one-dimensional Euler equations for compressible fluid flow \cite{wylie1978fluid}, which is a system of coupled partial differential equations (PDEs). That work was done in light of seminal and contemporary studies on the significance of the various terms in the Euler equations under the physical regime of normal pipeline flow \cite{Osiadacz1984,misra2020monotonicity}.  Direct simulation of transient flow PDEs in pipelines on a scale of thousands of kilometers is computationally challenging, and hence developing reduced-order ordinary differential equation (ODE) or differential algebraic equation (DAE) models for these systems is an active area of research \cite{Chapman2005}.  Recent studies developed a reduced order DAE model for networks of horizontal pipes \cite{Grundel2013}, and the extension to optimal control of intraday operations by gas compressors \cite{Zlotnik2015a}. Since then, similar models have been applied for joint power market optimization and pipeline management \cite{Zlotnik2016,Roald2020}, state and parameter estimation \cite{Sundar2018,Sundar2019}, and intra-day congestion pricing in gas pipeline systems \cite{Zlotnik2019}.

Previous studies have sought to apply linearization to the gas flow PDEs about the steady-state mass flow rates and pressures, in order to obtain linear state space models for gas pipeline dynamics \cite{Reddy2006,Alamian2012}. The motivation was to apply traditional techniques for state estimation and control that are available for linear systems. However, as noted above, the emerging influence of gas-fired power generation causes a wide range of transient phenomena in which flows throughout the supplying pipeline networks deviate substantially from the behavior approximated by steady-state models. Moreover, the approaches are generally not tested for computational scalability, and the only components for which reduced order models have been built include horizontal pipes and compressors.  To the best of our knowledge, \emph{there are no published studies that include physics-based spatio-temporal models of underground gas storage facilities into the operational transient pipeline optimization problem}.

Natural gas storage facilities have been included into the pipeline operational optimization problem in the form of simple models of storage as a source or a sink with capacity limits \cite{He2017}, i.e., gas can be withdrawn from or injected into the the storage at any rate as long as it is within its capacity limits. Such models don't resolve the complex spatio-temporal physics of large subsurface formations, such as depleted petroleum reservoirs, that are typically exploited by these facilities.  On the other hand, there are studies that model subsurface physics of single reservoirs at high fidelity \cite{Basumatary2005,Oldenburg2001}. The latter studies examine the complex interactions between the reservoir and the process of injecting and withdrawing gas from it. The complexity of these models makes their direct adaptation into the intra-day operational problem for pipeline management intractable.  Thus, models for storage are either simplistic \cite{He2017}, treating the storage as a source term with specified output and thus completely neglecting the physics, or are too complicated to be useful in the present setting. Therefore, there is a compelling need to develop reduced order models for a larger set of components in  gas pipeline systems that model the transient dynamics in a way that is amenable to integration with the associated optimal control problems.

In this study, we develop reduced order DAE models for all the components in a natural gas pipeline network including models for rectilinear, non-horizontal pipes and storage facilities. With regards to storage facilities, there are predominantly three types: depleted oil and gas reservoirs, salt caverns, and aquifers,  with each type having its own advantages and disadvantages. Around 78\% of the storage facilities in the US are depleted oil and gas reservoirs \cite{Lokhorst2005}, and thus we focus on developing a reduced-order model for this type of storage only. \edit{We provide a brief review of underground gas storage modeling and recent studies on optimization and control of reduced-order models of such systems later when we present the storage models in Sec. \ref{subsec:storage}.} The DAE model for a rectilinear, non-horizontal pipe is another novel modeling component of this paper, i.e, we generalize the available reduced-order models for horizontal pipes \cite{Sundar2018}. This generalized model is put to effective use in characterizing the physics of the injection and withdrawal wells used in the storage facility. The physical models of each component in the pipeline system are then put together and a nonlinear optimal control problem for performing day-to-day operations is formulated for a general gas pipeline network. The objective function of the optimal control problem is chosen to represent both economic and operational factors. The formulation is then solved to local optimality using standard gradient descent techniques. The efficacy of the formulations and algorithms are demonstrated using extensive simulation and computational experiments on a 6-junction test network and a model based on an actual pipeline network in the United States.

The rest of the paper is organized as follows: In Section \ref{sec:modeling}, we describe the newly developed reduced order model for gas storage facilities, and its synthesis with gas pipeline system models.  Section \ref{sec:nlocp} describes the optimal control problem of flow scheduling over a pipeline system with underground storage.  The developed problem formulation is applied to several test systems in a computational study described in Section \ref{sec:results}, and we conclude in Section \ref{sec:conclusion}.

\section{Modeling} \label{sec:modeling}

Here we extend the previously developed gas pipeline models used for transient optimization \cite{Sundar2018,Zlotnik2015} to include subsurface storage facilities and non-horizontal pipes.

\subsection{Gas Pipeline Dynamics} \label{subsec:pde}
The flow of compressible natural gas within a pipeline under isothermal conditions is described by specifying equations that govern i) the evolution of field quantities that represent the pressure $p(x, t)$ (or the density $\rho(x, t)$) and the velocity $v(x, t)$ of the gas as a function of time $t \in [0, T]$ and the spatial location $x \in [0, L]$ within the pipeline, and ii) the equation of state that relates pressure and density through the compressibility of the gas.  We fix an origin at one end of the pipe and designate the direction along the pipe to be the positive coordinate axis, so that  the velocity  is  also assumed to be directed along it.  
With this convention, we designate the mass flux $\varphi:= \rho v$ (units \si{\kilogram\per\square\metre\per\second}), and the mass flow at a location as $f := \varphi A$ (units \si{\kilogram\per\second}), where $A$ is the cross-sectional area available to the flow. The evolution of flow in a rectilinear pipeline is described by the system of one-dimensional Euler equations \cite{Thorley1987} 
\begin{subequations}
\begin{flalign}
& \partial_t \rho + \partial_x \varphi = 0, \label{eq:mass-balance} \\ 
& \partial_t \varphi + \partial_x (\rho v^2) = -\partial_x p + \rho g_{\parallel}  -\frac{\lambda}{2D} \frac{\varphi | \varphi |}{\rho}, \label{eq:momentum-balance}
\end{flalign}
\label{eq:euler}
\end{subequations}

\noindent
where the quantities $\lambda > 0, D > 0, g_{\parallel}$ are respectively the Darcy-Weisbach friction factor,  the diameter of the pipeline, and the component of the acceleration due to gravity \emph{along} the coordinate axis. Note that the value of $g_{\parallel}$ is determined by the orientation of the pipeline; for a pipe aligned at an angle $0 \leqslant \theta \leqslant \pi/2$ with the horizontal, we have $g_{\parallel} = \pm g\sin{\theta}$,  with the the sign determined by the choice of coordinate direction.
Equations \eqref{eq:mass-balance} and \eqref{eq:momentum-balance} represent the balance of  mass and momentum respectively, in  the regime when the changes in the boundary conditions are sufficiently slow so as to not form waves or shocks \cite{Thorley1987}. We direct the reader to a seminal work to justify the use of these equations for large-scale gas pipeline flows in the normal operating regime \cite{Osiadacz1984}, and a recent study on the friction-dominated approximation \cite{misra2020monotonicity}. We also assume that the density and pressure are related by the equation of state 
\begin{flalign}
p = ZRT \rho = a^2 \rho, \label{eq:eq-of-state}
\end{flalign}
where $Z$, $R$, $T$, and $a$ represent the gas compressibility factor, universal gas constant, (constant) temperature, and the (wave) speed of sound, respectively. 
While one could assume more complex, nonlinear relationships between the pressure and density, those are beyond the scope of the current work, as are  the incorporation of other effects caused by temperature, natural gas composition, and more realistic equation of state models \cite{Menon2005,Gyrya2019}. 

In Eq. \eqref{eq:momentum-balance}, the term $\partial_t \varphi$ is known to be much smaller than $\partial_x (\rho v^2) + \partial_x p$, and the ratio of $\partial_x p$ and $\partial_t \varphi$ is typically on the order of $1 : 10^{-3}$ \cite{Osiadacz1984}. Furthermore, the flow velocities are much smaller than the speed of sound, i.e., $v \ll a$. Hence, we may omit the terms $\partial_t \varphi, \partial_x (\rho v^2)$ from the equation. With these observations, the governing equations \eqref{eq:euler} simplify to 
\begin{subequations}
\begin{flalign}
& \partial_t \rho + \partial_x \varphi = 0, \label{eq:mass-balance-1} \\ 
& \partial_x p = \rho g_{\parallel}  -\frac{\lambda}{2D} \frac{\varphi | \varphi |}{\rho}, \label{eq:momentum-balance-1} \\ 
& p = a^2 \rho. \label{eq:eq-of-state-1}
\end{flalign}
\label{eq:euler-1}
\end{subequations}
There are many studies in the literature that validate and verify the use of Eq. \eqref{eq:momentum-balance-1} in place of Eq. \eqref{eq:momentum-balance}. An interested reader is referred to \cite{Osiadacz1984,Sundar2018,Dyachenko2017,Zlotnik2015,Mak2016} for these studies. Hence, throughout the rest of the article, gas flow dynamics through a pipeline are represented using Eqs. \eqref{eq:euler-1}. 
Note that while the system \eqref{eq:euler} is hyperbolic, the simplified system \eqref{eq:euler-1} is parabolic \cite{Sundar2018}, and  is well-posed when the initial conditions are supplemented by boundary conditions, consisting of prescribed values for one of $\rho(0, t)$ or $\varphi(0, t)$ and one of $\rho(L, t)$ or $\varphi(L, t)$. Gas pipeline managers typically aim to maintain the system state in synchronization with loads, which are roughly periodic at $24$ hours. The time periodicity dictates our modelling of other components of a full pipeline network. We examine appropriate boundary and endpoint conditions for optimal control of a pipeline network in Section~\ref{sec:nlocp}.

We non-dimensionalize the system of equations \eqref{eq:euler-1} in order to improve numerical conditioning by applying the transformations
\begin{flalign}
\hat t = \frac{t}{t_0}, ~ \hat x = \frac{x}{\ell},~ \hat p = \frac{p}{p_0}, ~\hat{\varphi} = \frac{\varphi}{\varphi_0}, \label{eq:non-dimensional-constants}
\end{flalign}
where $p_0$ and $\ell$ are respectively the nominal length and pressure, and the relationship between $p_0$, $\ell$, $\varphi_0$, and $t_0$ is given by $\rho_0 = p_0/a^2$, $\varphi_0 = a\rho_0$, and $t_0 = \ell/a$. The resulting system of equations after non-dimensionalization is given by 
\begin{flalign}
\frac{\partial \hat{\varphi}}{\partial \hat x} + \frac{\partial \hat{\rho}}{\partial \hat t} = 0,~ 
\frac{\partial \hat{p}}{\partial \hat x} = \frac{\ell}{a^2}\hat{\rho} g_{\parallel}  -\frac{\lambda \ell}{2D} \frac{\hat{\varphi}|\hat{\varphi}|}{\hat{\rho}},~\hat{p}  = \hat{\rho}. \label{eq:euler-nd-hat}
\end{flalign}
For the sake of readability, we henceforth omit the hat symbols from the non-dimensionalized  system of equations and write 
\begin{flalign}
\partial_x \varphi+ \partial_t \rho = 0 \text{ and } \partial_x \rho^2 = 
\frac{2\ell g_{\parallel}}{a^2} \rho^2 -\frac{\lambda \ell}{D} {\varphi|\varphi|}. \label{eq:euler-nd}
\end{flalign}

\subsection{Gas Pipeline Network - Notations} \label{subsec:notations} 
A gas transmission pipeline network may be modelled as a graph that consists of pipes and compressors (edges) interconnected at junctions (nodes) where gas can be injected into the system (intake points), withdrawn from the system (off-take points), or interact with storage facilities.  The operational state of a storage facility enables gas to be injected or withdrawn from the facility, and this is determined based on seasonal needs. We designate the gas transmission pipeline network as a connected directed graph $\mathcal G = (\mathcal V, \mathcal P \cup\, \mathcal C)$ where $\mathcal V$ represents the set of junctions in the pipeline network and the sets $\mathcal P$ and $\mathcal C$ represent the respective set of horizontal pipelines and the set of compressors that connect two junctions. 
It is reasonable to assume that the values of gas density at pipeline junctions are known at all times, thus
each junction $i \in \mathcal V$ is associated with a time-varying density variable $\rho_i(t)$. Also associated with each junction $i \in \mathcal V$ are three sets $\mathcal R_i$, $\mathcal  D_i$, and $\mathcal S_i$, which denote the set of intake points, the set of off-take points, and the set of storage units attached to the junction $i$. We set $\mathcal R = \cup_i \mathcal R_i$, $\mathcal D = \cup_i \mathcal D_i$, and $\mathcal S = \cup_i \mathcal S_i$. 

Each  pipeline $(i, j) \in \mathcal P$ that connects two junctions $i$ and $j$ is associated with four parameters $L_{ij}$, $D_{ij}$, $\lambda_{ij}$, and $\theta_{ij}$ that denote the pipe's length, diameter, friction factor, and orientation angle respectively. For an individual pipe, it is conceivable that mass flux measurements are available  at its two ends, say $\varphi_{ij}^{in}(t)$ and $\varphi_{ij}^{out}(t)$. Based on these two measurements, one may define two quantities for each pipe $(i, j) \in \mathcal P$,
\begin{flalign}
\varphi^{+}_{ij}(t) = \frac{\varphi^{out}_{ij}(t) + \varphi^{in}_{ij}(t)}{2}, ~ \varphi^{-}_{ij}(t) = \frac{\varphi^{out}_{ij}(t) - \varphi_{ij}^{in}(t)}{2}, \label{eq:avg-def}
\end{flalign} 
that will be convenient for later use.

We let $\rho_{ij}(x, t)$ and $\varphi_{ij}(x, t)$ denote the density and mass flux throughout the pipeline.
As described earlier in the context of Eq. \eqref{eq:euler}, for each pipeline $(i,j)$, we choose the  coordinate direction $i\rightarrow j$, and  the sign of the the mass flux $\varphi_{ij}(x, t)$ then indicates whether flow is along or against the assumed direction. It naturally follows that
\begin{subequations}
\begin{flalign}
&\rho_{ij}(0, t) = \rho_i(t),~~ \rho_{ij}(L_{ij}, t)  = \rho_j(t), \label{eq:rho-bdry} \\ 
&\varphi_{ij}(0, t) = \varphi_{ij}^{in} = \varphi^{+}_{ij}(t) - \varphi^{-}_{ij}(t), \text{ and } \label{eq:in} \\ 
&\varphi_{ij}(L_{ij}, t) = \varphi_{ij}^{out} = \varphi^{+}_{ij}(t) + \varphi^{-}_{ij}(t). \label{eq:out}
\end{flalign}
\label{eq:bdry}
\end{subequations}
As gas flows through a long pipeline, the friction in the pipeline causes pressure to gradually decrease along the direction of flow. Hence, the network consists of compressors at regular intervals of pipeline length that boost the line pressure to meet the minimum pressure requirements for delivery of gas. These compressors are controllable actuators that connect two junctions; furthermore, we assume a compressor's length is negligible since the size of a compressor is small relative to the length of a  pipeline. Hence, the flow through a compressor $(i, j) \in \mathcal C$ that connects junction $i$ to $j$ is modeled as a multiplicative increase in density of the gas as it flows from  junction $i$ to junction $j$, i.e., $\rho_j(t) = \alpha_{ij}(t) \cdot \rho_i(t)$. Here, $\alpha_{ij}(t)$ represents the time-dependent compression ratio between the suction density, $\rho_i(t)$, and discharge density, $\rho_j(t)$. 

In the subsequent sections, we derive a reduced-order model that includes each of the above network components based on the system of equations \eqref{eq:euler-nd}. We first start by deriving such a model for short segments of rectilinear pipes, leveraging previous work  \cite{Grundel2013,Zlotnik2015a,Sundar2018}. 

\subsection{Rectilinear Pipe} \label{subsec:short_pipes}
Rectilinear pipes in the gas transmission network are used to transport gas from one location to another. To derive a reduced-order model for a single pipe in the transmission network that connects two junctions, we first observe that for a sufficiently small length of the pipeline segment, the relative density and mass flux difference between the neighboring junctions is small at all times, i.e., for a small enough pipeline $(i, j)$ of length $L_{ij} \leqslant \Delta$, we have
\begin{flalign}
\left|\frac{\rho_j(t) - \rho_i(t)}{\rho_j(t) + \rho_i(t)}\right| \ll 1, \left|\frac{\varphi^{-}_{ij}(t)}{\varphi^{+}_{ij}(t)} \right| \ll 1, ~ \forall t \in [0, T]. \label{eq:delta}
\end{flalign}
Hence, without loss of generality, we assume that any  pipeline in the transmission system is of length at most $\Delta$ units (if the pipe is actually longer than $\Delta$, we segment the pipe into a set of sub-pipes by adding junctions such that each sub-pipe has length at most $\Delta$). 
 
Hence, the system of equations that govern the flow of natural gas in the pipe is given by 
\begin{subequations}
\begin{flalign}
\partial_x \varphi_{ij} + \partial_t \rho_{ij} = 0, \label{eq:mass-r-pde} \\ 
\partial_x \rho_{ij}^2  - \frac{2\ell g_{\parallel}}{a^2}\rho_{ij}^2 = -\frac{\lambda_{ij} \ell}{D_{ij}} \varphi_{ij}|\varphi_{ij}|. \label{eq:momentum-r-pde}
\end{flalign}
\label{eq:pipe-pde}
\end{subequations}
To build a reduced-order model, we first integrate the conservation equations in Eq. \eqref{eq:pipe-pde} along the length coordinate:
\begin{subequations}
\begin{align}
\!\! & \int_0^{{L}_{ij}} \left(\partial_x \varphi_{ij} + \partial_t \rho_{ij}\right) dx = 0,  \label{eq:mass-int} \\ 
\!\! & \int_0^{{L}_{ij}} \!\!\left( \partial_x \rho_{ij}^2 \!-\! \frac{2\ell g_{\parallel}}{a^2}\rho_{ij}^2 \right)\!dx \! = \! -\frac{\lambda_{ij} \ell}{D_{ij}} \int_0^{{L}_{ij}} \!\!\varphi_{ij}|\varphi_{ij}| dx. \label{eq:momentum-int}
\end{align}
\label{eq:integrate-pde}
\end{subequations}
Wherever possible, we integrate in closed form, else we use the mid-point rule of integration, the fundamental theorem of calculus, and averaging to evaluate the integrals in terms of the lumped quantities $\rho_i(t)$, $\rho_j(t)$, $\varphi^{+}_{ij}(t)$, and $\varphi^{-}_{ij}(t)$ to reduce the system of equations \eqref{eq:integrate-pde} to
\begin{subequations}
\begin{align}
& {L_{ij}} \cdot (\dot{\rho}_j(t) + \dot{\rho}_i(t)) = - 4 \cdot \varphi^{-}_{ij}(t),  \label{eq:mass-pipe-1} \\ & e^{\beta}\rho_j^2(t) \!-\! \rho_i^2(t) \!=\! -\frac{\lambda_{ij} \ell L_{ij}}{D_{ij}} \left(\!\frac{e^{\beta} -1}{\beta} \!\right) \cdot \varphi^{+}_{ij}(t) | \varphi^{+}_{ij}(t)|,  \label{eq:momentum-pipe-1}
\end{align}
\label{eq:pipe-ode-1}
\end{subequations}
where $\beta \triangleq -2\ell g_{\parallel} L_{ij}/a^2$, and the dot over a quantity denotes the time derivative.
The system of equations \eqref{eq:pipe-ode-1} is the reduced-order ordinary differential equation (ODE) model for \emph{any}  pipe $(i, j) \in \mathcal P$ in the network with length at most $\Delta$ units. 
For a horizontal pipeline, we have $g_{\parallel} = 0$, so that the limit $\beta \rightarrow 0$ in Eq. \eqref{eq:pipe-ode-1} yields
\begin{subequations}
\begin{align}
& {L_{ij}} \cdot (\dot{\rho}_j(t) + \dot{\rho}_i(t)) = - 4 \cdot \varphi^{-}_{ij}(t),  \label{eq:mass-pipe-horz} \\ 
& \rho_j^2(t)- \rho_i^2(t) = -\frac{\lambda_{ij} \ell L_{ij}}{D_{ij}} \cdot \varphi^{+}_{ij}(t) | \varphi^{+}_{ij}(t)|.  \label{eq:momentum-pipe-horz}
\end{align}
\label{eq:pipe-ode-horz}
\end{subequations}

\begin{figure}[htbp]
    \centering
    \includegraphics[scale=0.8]{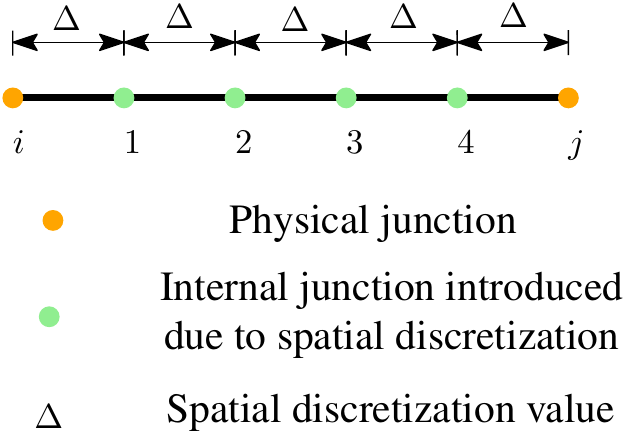}
    \caption{Illustration of the spatial discretization scheme for a rectilinear pipe.}
    \label{fig:pipe-discretization}
\end{figure}

\edit{In the computational implementations, each pipe is discretized along its axial direction. The differential equation models for rectilinear pipes in Eq. \eqref{eq:pipe-ode-1} are valid for a small enough pipeline $(i,j)$ with length $L_{ij} \leqslant \Delta$ as in Eq. \eqref{eq:delta}. The pipelines are modeled as one-dimensional objects and the meshing is applied as a partitioning of the axial dimension of each pipeline into segments (see Fig. \ref{fig:pipe-discretization}). This is a standard discretization approach in pipeline simulation \cite{Thorley1987}. Empirical studies that examine discretization lengths required for acceptable accuracy for pipeline flow modeling have been done previously. In particular, in \cite{Grundel2013} the authors  found that $10$ \si{\km}, i.e., $\Delta = 10$ \si{\km} is a sufficiently fine discretization for a pipeline of arbitrary length. Furthermore, point wise convergence for such a discretization scheme to the conservation equations in in Eq. \eqref{eq:pipe-pde} has been proven in previous work \cite{Zlotnik2015a}.}

% \vspace{-5ex}
\subsection{Compressor} \label{subsec:compressors}
Each compressor $(i, j) \in \mathcal C$ that connects junctions $i$ and $j$ is modeled as a device that performs a time-dependent multiplicative increase of density from junction $i$ to junction $j$. Thus the convention for flow direction through a compressor is not arbitrary, but part of the definition of the device. The multiplicative factor between the outlet pressure and inlet pressure, called the compression ratio, is denoted by $\alpha_{ij}(t) \geqslant 1$. Also, we let $f_{ij}(t)$ denote the mass flow rate through a compressor that is assumed to be constant throughout the length of the compressor; when $f_{ij}(t)$ is positive (negative), then the gas flows from $i$ to $j$ ($j$ to $i$). Given these possibilities, compressors can be of two predominant types: (i) compressor $(i, j)$ that allows flow in only one direction and boosts density ($i \rightarrow j$), i.e., $f_{ij}(t) \geqslant 0$ and (ii) compressor $(i, j)$ that allows bi-directional flow but boosts density in one direction ($i \rightarrow j$) only, i.e., $f_{ij}$ can be either positive or negative. 
For the compressor $(i, j) \in \mathcal C$ that allows flow in one direction only, the model is:
\begin{flalign}
\rho_j(t) = \alpha_{ij}(t) \cdot \rho_i(t) \text{ and } f_{ij}(t) \geqslant 0. \label{eq:c-1}
\end{flalign}
For the compressor $(i, j) \in \mathcal C$ that allows uncompressed reverse flow, the model is:
\begin{flalign}
\rho_j(t) = \alpha_{ij}(t) \cdot \rho_i(t) \text{ and } f_{ij}(t) \cdot (1 - \alpha_{ij}(t))  \leqslant 0. \label{eq:c-2}
\end{flalign}
In Eq. \eqref{eq:c-2}, when $f_{ij}(t) < 0$, then flow is from junction $j$ to $i$; in this case, the constraints in Eq. \eqref{eq:c-2} enforce $\alpha_{ij}(t) = 1$, i.e., gas is not compressed  when flow is in the reverse direction. 

We assume that the gas compression process is adiabatic, i.e., no heat transfer occurs between the gas and the surroundings. On this basis, the work done by a compressor $(i, j) \in \mathcal C$ in compressing $1$ kg of natural gas is given by 
\begin{flalign}
W_{ij}(t) = \left(\frac{\gamma \cdot T}{\gamma -1}\right) \left(\frac{286.76}{G}\right) \left( \alpha_{ij}(t)^{m} - 1 \right) \text{ \si{\joule\per\kilogram}}\label{eq:work-compressor}
\end{flalign}
where, $T$ is the temperature at which the gas is compressed, $G$ is gas gravity (dimensionless), and $\gamma$ is the ratio of specific heats of natural gas, and $m = (\gamma-1)/\gamma$ \cite{Menon2005}. Given Eq. \eqref{eq:work-compressor}, the power consumed by the compressor $(i, j) \in \mathcal C$ to compress $f_{ij}(t)$ units of flow in \si{\kilogram\per\second} at a ratio $\alpha_{ij}(t)$ is given by $W_{ij}(t) \cdot f_{ij}(t)$ or $W_{ij}(t) \cdot f_{ij}(t) \cdot \varphi_0$, when $f_{ij}(t)$ is non-dimensionalized by nominal mass flux and unit area, i.e,  by $\varphi_0 \cdot 1 \cdot 1$. 

\subsection{Storage Facility} \label{subsec:storage} 
A natural gas storage facility typically consists of several components such as: i) underground reservoirs, ii) injection and withdrawal wells, iii) compressors and regulators, iv) gathering systems, v) dehydrators, vi) metering facilities, etc. For planning midstream pipeline operations, one is interested in the amount of gas that can be stored, and the rate at which the gas can either be injected into or withdrawn from the storage. Hence the components of a storage facility that we model are: i) reservoir (where gas is stored); ii) wells (used to inject or withdraw gas from the reservoir), and iii) a compressor and a regulator station (to control the pressure at the well head).

The amount of gas that can be stored in a reservoir depends  primarily on the volume of the reservoir, but it is also affected by other physical phenomena, such as the deformability of the boundary walls, or physical/chemical absorption. 
Similarly, the rate at which gas can be withdrawn from or injected into the reservoir depends on several factors such as the pressure gradient across the well and the reservoir, friction along the wells, permeability of the reservoir to gas, viscosity of gas, etc. However, given the context of this article, a high-fidelity model accounting for \emph{all} these factors is not viable for the following two reasons: i) lack of readily available data on reservoir characteristics, and ii) the computational intractability of the transient optimization problems that would result from a complicated model. Hence, in this article, we consider a simplified model for the reservoir and the storage as a whole that does not compromise significantly on the physics, yet is tractable for computational optimal control.

\edit{For a comprehensive review of the modeling, analysis, and engineering of underground gas storage facilities that has been used for the design, implementation, and operation of these sites in the 20th century, we refer the reader to the authoritative work of Orin Flanigan \cite{flanigan1995underground}. These practical approaches were developed before the widespread use in the oil and gas industry of complex methods such as machine learning and multi-phase flow modeling, and almost all current subsurface gas storage facilities were built with and are operated based on simple engineering models that rely on empirically derived relationships. More recent scientific studies have sought to develop reduced-order models and apply optimization and optimal control to complex subsurface structures \cite{jansen2008model,chen2020frankenstein}.  In general, advanced problems such as optimal control are very challenging to tackle using sophisticated full-physics simulation models because they require many parameters, and are difficult to express in ways that are amenable to optimization solvers.  Therefore, we use the historically proven approach of simplified modeling in order to make a first step in integrating subsurface storage models into transient optimization of pipeline systems, expecting that future studies will integrate more advanced models and validation case studies.
}

\begin{figure}[htbp]
    \centering
    \includegraphics[scale=0.85]{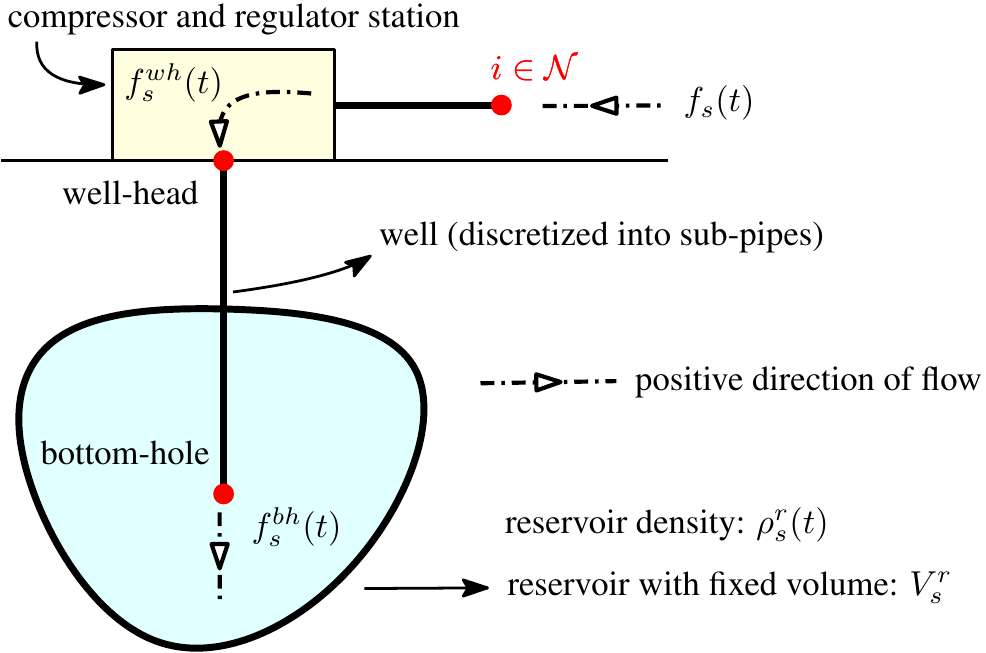}
    \caption{Storage schematic. The storage unit is connected to the junction $i \in \mathcal N$. Note that positive direction of flow indicates injection into reservoir. In every case, $f_s(t) = f_s^{wh}(t)$.}
    \label{fig:storage}
\end{figure}

A schematic of the components considered in a storage system is shown in Fig. \ref{fig:storage}. Given a storage system $s \in \mathcal S$ attached to a junction $i \in \mathcal V$, i.e., $s \in \mathcal S_i$, we first associate a flow variable $f_{s}(t)$. This variable denotes the amount of  gas (specified by the mass per unit time) exchanged between the pipeline network at junction $i$ and the storage system. Depending on the sign, i.e., if $f_s(t) > 0$, the gas is injected into the storage, else if $f_s(t)<0$, it is withdrawn from the storage. Next, we describe the models for the reservoir and the well. 

\subsubsection{Reservoir and Wells} \label{subsubsec:res-well}
A simplified storage model consists of an underground reservoir that is used to store the gas and a vertical pipe or a well through which gas can be injected into or withdrawn from the reservoir. In reality, a single storage unit has multiple dedicated injection and withdrawal wells, but for ease of modeling, we use a single, dual-purpose well, since injection and withdrawal do not happen simultaneously in practice. We make the following modeling assumptions about the reservoir: i) The reservoir is modelled as a domain with fixed boundary geometry, so that it has constant volume, denoted by $V_s^r$; ii) the gas is distributed homogeneously within the reservoir with a time-varying density denoted by $\rho_s^r(t)$ that increases (decreases) as gas is injected into (withdrawn from) the reservoir; and iii) the injection/withdrawal processes cause negligible temperature changes, so that the gas is stored at constant temperature. 

The junctions at the top and the bottom of the well are referred to as the well-head and the bottom-hole, respectively. 
The well, being a vertical pipe, has a diameter, length, and friction-factor parameter associated with it and gas flows through the well due to the net effect of gravity and the pressure gradient driving the flow. Each well is discretized into a finite number of sub-pipes such that each sub-pipe is of length at most $\Delta$ units. For all sub-pipes in the well, we assume that the coordinate axis is vertically downward. Given these assumptions, the governing equations for any sub-pipe in the well $(i, j)$, are given by Eq. \eqref{eq:pipe-ode-1} with $g_{\parallel} = g$, i.e., $\beta = \frac{2\ell gL_{ij}}{a^2}$. 

The amount of gas injected into or withdrawn at the bottom-hole and well-head is denoted by $f_s^{bh}(t)$ and $f_s^{wh}(t)$, respectively. With our established convention,  $f_s^{bh}(t) > 0$ and $f_s^{wh}(t) > 0$ when gas is being injected into the reservoir and negative otherwise. The density at the bottom-hole is assumed to be the same as the density of gas in the reservoir, $\rho_s^r(t)$. Then, the following equation governs the change in reservoir density as a function of the injection or withdrawal:
\begin{flalign}
V_s^r \cdot \dot{\rho}_s^r(t) = f_s^{bh}(t). \label{eq:reservoir}
\end{flalign}
The Eq. \eqref{eq:reservoir} is in non-dimensional form when the non-dimensional constant for volume is given by $V_0 = \ell \cdot 1 \cdot 1$.

To maintain minimum density required for withdrawal, a certain amount of gas must always be present in the reservoir. This gas is referred to as base gas and is not allowed to be withdrawn. Similarly, to maintain the integrity of the reservoir, the maximum amount of gas that can be stored is also restricted. If $m_s^{\min}$ and $m_s^{\max}$ are the minimum and maximum amounts of gas that can be stored in the reservoir of storage $s \in \mathcal S$, then the mass in the reservoir at time instant $t$ must satisfy the following non-dimensional inequality constraint:
\begin{flalign}
m_s^{\min} \leqslant \rho_s^r(t) \cdot V_s^r \leqslant m_s^{\max}. \label{eq:reservoir-limits}
\end{flalign}
Here, the non-dimensional constant for mass $m_0 = \rho_0 \cdot V_0$. 

\subsubsection{Compressor and Regulator Station} \label{subsubsec:cr-station}
Often, the well-head pressure required to obtain a specified amount of gas flow rate is beyond acceptable pressure levels at the junction to which the well is connected (junction $i \in \mathcal N$ in Fig. \ref{fig:storage}). Therefore, there is a need to either increase or decrease the well-head pressure during injection and withdrawal phases, respectively. In practice, this is achieved through separate compressor stations and density-reduction (regulator) stations. For ease of modeling, we combined these two components into a single component that can perform both compression or density-reduction. This is achieved by introducing a multiplicative factor, $\alpha_s(t)$ that relates the density at the well-head  ($\rho_s^{wh}(t)$) with the density $\rho_i(t)$ at the junction to which the storage is attached to, according to the following relation
\begin{subequations}
\begin{flalign}
\rho_i(t) = \alpha_s(t) \cdot \rho_s^{wh}(t), \label{eq:boost} \\ 
\frac{1}{\alpha_s^{\max}} \leqslant \alpha_s(t) \leqslant \alpha_s^{\max} \label{eq:limits}
\end{flalign}
\label{eq:cr-station}
\end{subequations}
where, $\alpha_s^{\max} > 1$ is the maximum compression that can be achieved at the compressor station. The maximum density-reduction that is achieved is given by $\frac{1}{\alpha_s^{\max}}$. The flow through the compressor and regulator station is given by $f_s^{wh}(t)$, which is also equal to the flow coming out of the junction $i$, i.e., $f_s(t)$. 

\subsection{Gas Intake and Off-take Points} \label{subsec:rd} 
We now define the variables that are associated with the gas intake and off-take points, which are located at junctions and at which gas is injected into and/or withdrawn from the pipeline system. A flow variable $f_r(t) > 0$ ($f_d(t)>0$) is associated with each intake (off-take) point $r \in \mathcal R$ ($d \in \mathcal D$), and represents the amount of gas that can be injected into (withdrawn from) the pipeline system. The only constraint on these flows are the injection and withdrawal limits given by $f_r^{\max}(t)$ and $f_d^{\max}(t)$, respectively. 

Having outlined the governing equations for each component in a pipeline network, in the next section we synthesize them within an optimal control problem. 

\section{Nonlinear Optimal Control Problem} \label{sec:nlocp} 
To formulate the nonlinear optimal control problem for flow scheduling on a gas pipeline network with storage, we first introduce a few definitions. Recall that for the construction of the reduced-order model, we discretize long pipes into sub-pipes of maximum length $\Delta$. Thus, the pipeline network is represented as a directed graph $\mathcal G = (\mathcal V, \mathcal E)$ with a length $L_{ij}$ associated with each edge $(i,j) \in \mathcal E$ such that $L_{ij} \leq \Delta$. Here, $\mathcal V$ is the set of junctions in the pipeline network, and the set  of  pipelines $\mathcal P$  and  compressors $\mathcal C$  that connect  two  junctions comprise  $\mathcal E$. We assume that all pipelines in $\mathcal P$ are horizontal.

We now append  a set $\mathcal V^s$ that includes  junctions comprising each well of each storage unit to the original set of junctions in $\mathcal V$. Similarly the set of pipes and compressors in $\mathcal P$ and $\mathcal C$ is appended with the set of wells $\mathcal W$, one for each storage, and the set of compressor and regulator stations $\mathcal C^s$. The augmented gas pipeline network is $\mathcal G^{\dagger} = (\mathcal V \cup \mathcal V^s, \mathcal P \cup \mathcal C \cup \mathcal W \cup \mathcal C^s)$, and the optimal control problem is formulated for ${\mathcal G}^{\dagger}$ subsequently. 

\subsection{Decision Variables} \label{subsec:variables}
The decision variables to be optimized are as follows:
\begin{subequations}
\begin{flalign}
\forall i \in {\mathcal V} \cup {\mathcal V}^s: \quad & \rho_i(t) \label{eq:junction_var} \\
\forall (i, j) \in {\mathcal P} \cup {\mathcal W}: \quad  & \varphi_{ij}^+(t), \, \varphi_{ij}^-(t), \, \varphi_{ij}^{in}(t), \, \varphi_{ij}^{out}(t) \label{eq:pipe_var} \\ 
\forall (i, j) \in \mathcal C: \quad & \alpha_{ij}(t), \, f_{ij}(t), \, W_{ij}(t) \label{eq:compressor_var} \\ 
\forall (i, j) \in \mathcal C^s: \quad & \alpha_{ij}(t) \label{eq:cr_var} \\ 
\forall r \in \mathcal R: \quad & f_r(t) \label{eq:receipt_var} \\ 
\forall d \in \mathcal D: \quad & f_d(t) \label{eq:delivery_var} \\ 
\forall s \in \mathcal S: \quad & \rho_s^r(t), \, f_s^{bh}(t), \, f_s^{wh}(t), \, f_s(t) \label{eq:storage_var}
\end{flalign}
\label{eq:variables}
\end{subequations}
The variables listed in  Eq. \eqref{eq:junction_var} - Eq. \eqref{eq:pipe_var} include density variables for each junction in \eqref{eq:junction_var} and flux variables associated for each pipe in the network in \eqref{eq:pipe_var}.
The variables associated with a compressor, i.e., the multiplicative boost, the mass flow, and the work done in compression of  $1$ \si{\kilogram} of gas are listed in  Eq. \eqref{eq:compressor_var}.
Finally, Eq. \eqref{eq:cr_var} lists the variables associated with a compressor and regulator station at a storage, while  Eq. \eqref{eq:receipt_var} -- \eqref{eq:storage_var} enumerate the variables used for gas intake points, off-take points, and storage respectively. 

\subsection{Constraints} \label{subsec:constraints} 
We classify the constraints of the problem into i) operational constraints and ii) component physics constraints. The operational constraints are given by 
\begin{subequations}
\begin{flalign}
\forall i \in \hat {\mathcal V} \subset \mathcal V: \quad & \rho_i(t) \text{ is specified} \label{eq:slack_rho} \\
\forall i \in {\mathcal V}: \quad & \rho^{\min} \leqslant \rho_i(t) \leqslant \rho^{\max} \label{eq:rho_junction_limits} \\
\forall i \in {\mathcal V}^s: \quad & \rho_s^{\min} \leqslant \rho_i(t) \leqslant \rho_s^{\max} \label{eq:rho_well_junction_limits} \\
\forall (i, j) \in {\mathcal P}: \quad &  -\varphi^{\max} \leqslant \varphi_{ij}^{in}(t), \varphi_{ij}^{out}(t) \leqslant \varphi^{\max}\label{eq:pipe_flux_limits} \\
\forall (i, j) \in \mathcal C: \quad & -f_{ij}^{\max} \leqslant f_{ij}(t) \leqslant f_{ij}^{\max}  \label{eq:compressor_flow_limits} \\ 
\forall (i, j) \in \mathcal C: \quad & 1 \leqslant \alpha_{ij}(t) \leqslant \alpha^{\max}  \label{eq:compressor_ratio_limits} \\ 
\forall (i, j) \in \mathcal C: \quad & W_{ij}(t) \cdot f_{ij}(t) \cdot \varphi_0 \leqslant P^{\max} \label{eq:compressor_power_limits} \\ 
\forall r \in \mathcal R: \quad & 0 \leqslant f_r(t) \leqslant f_r^{\max}(t) \label{eq:receipt_limits} \\ 
\forall d \in \mathcal D: \quad & 0 \leqslant f_d(t) \leqslant f_d^{\max}(t) \label{eq:delivery_limits} \\ 
\forall s \in \mathcal S: \quad & m_s^{\min} \leqslant \rho_s^r(t) \cdot V_s^r \leqslant m_s^{\max} \label{eq:storage_limits} \\ 
\forall s \in \mathcal S: \quad & -f_s^{\max} \leqslant f_s(t) \leqslant f_s^{\max} \label{eq:storage_flow_limits} \\
\forall (i, j) \in \mathcal C^s: \quad & \frac{1}{\alpha_s^{\max}} \leqslant \alpha_s(t) \leqslant \alpha_s^{\max} \label{eq:cr_limits}
\end{flalign}
\label{eq:operation-limits}
\end{subequations}
The first constraint in Eq. \eqref{eq:slack_rho} fixes the density at certain nodes in the system that are referred to as slack nodes, which model intake locations where densities are fixed externally.
For the remaining inequality constraints, the bound values are specified as a part of the input data along with the gas pipeline network model. The constraints \eqref{eq:rho_junction_limits} and \eqref{eq:rho_well_junction_limits} delimit densities at pipeline junctions and well junctions, respectively. The constraints \eqref{eq:pipe_flux_limits} and \eqref{eq:compressor_flow_limits} delimit flows through pipes (at endpoints) and through compressor stations.  Constraints \eqref{eq:compressor_ratio_limits} and \eqref{eq:compressor_power_limits} limit the density boost ratio and applied power for the compressors, where $P^{\max}$ denotes the maximum compressor operating power. The constraints \eqref{eq:receipt_limits} and \eqref{eq:delivery_limits} limit the nodal receipts and deliveries of gas throughout the pipeline network, where the constraint bounds $f_r^{\max}(t)$ and $f_d^{\max}(t)$ are defined by the time-dependent receipt and delivery nominations provided by users of the pipeline network.  Finally, the constraints \eqref{eq:storage_limits}, \eqref{eq:storage_flow_limits}, and \eqref{eq:cr_limits} delimit the quantity of gas in storage, as well as the operating range of flow rate and wellhead compression ratio for injection or withdrawal actions.

The physics constraints divide naturally into three groups: boundary conditions,  nodal balance constraints (balance of mass at junctions), and local edge constraints (governing equations for pipes, wells or compressors), and we shall outline each of them in turn. \\

\noindent \textit{Endpoint conditions}: In this study, we use a time horizon of $T = 24$ hours for the optimal control problem. Gas pipeline operators require that the system, starting at $t=0$, return to a nominal state at $t=24$ such that the line-pack (the total mass of gas) values within local subsystems of the pipeline are restored. In our implementation, this requirement is met by enforcing time-periodic boundary conditions on a subset of state variables. Specifically, time-periodicity is enforced on the variables  $\rho_i(t)$ in Eq. \eqref{eq:junction_var}, $\varphi_{ij}^{in}(t)$ and $\varphi_{ij}^{out}(t)$ in Eq. \eqref{eq:pipe_var}, and 
$f_{ij}(t)$ in Eq. \eqref{eq:compressor_var} by constraining these functions to take the same value at $t = 0$ and $T=24$ hours, according to
    \begin{align} \label{eq:endpoint_cond}
        \rho_i(0) &= \rho_i(T), \,\, \varphi_{ij}^{in}(0) = \varphi_{ij}^{in}(T), \,\, f_{ij}(0) = f_{ij}(T).
    \end{align}
Because of their mutual relationship, all variables in Eq. \eqref{eq:pipe_var} and Eq. \eqref{eq:compressor_var} will satisfy the same periodicity condition.
 
\noindent \textit{Local edge constraints}: The governing equations of node connecting edge components have been discussed above individually. The dynamic (physics) constraints for pipes, wells, compressors, and wellhead compressors are
\begin{subequations}
\begin{flalign*}
& \forall (i, j) \in {\mathcal P}: \text{ Eq. \eqref{eq:pipe-ode-horz}, \eqref{eq:in}, \eqref{eq:out}, \eqref{eq:pipe-ode-horz},} &\\
& \forall (i, j) \in {\mathcal W}: \text{ Eq. \eqref{eq:in}, \eqref{eq:out}, Eq. \eqref{eq:pipe-ode-1} with $g_{\parallel} = g$,} &\\ 
& \forall (i, j) \in \mathcal C: \text{ Eq. \eqref{eq:c-1} or \eqref{eq:c-2}, Eq. \eqref{eq:work-compressor},} & \\ 
& \forall (i, j) \in \mathcal C^s: \text{ Eq. \eqref{eq:boost} }.
\end{flalign*}
\label{eq:physics-constraints}
\end{subequations}

\noindent \textit{Nodal balance constraints}:  Given a junction $i \in {\mathcal V} \cup {\mathcal V}^s$, we denote ``incoming'' pipes by $\delta_p^-(i):= \{j|(j,i) \in \mathcal{P} \cup \mathcal{W}$, and correspondingly, ``outgoing'' pipes by $\delta_p^+(i) := \{j|(i,j) \in \mathcal{P} \cup \mathcal{W}\}$. We define the sets $\delta_c^-(i)$, and $\delta_c^+(i)$ corresponding to the set of compressors $\mathcal{C}$ analogously.  Then, a nodal balance constraint for any junction $i \in {\mathcal V}$ can be stated as
\begin{flalign}
& \sum_{j \in \delta_p^+(i)}\varphi_{ij}^{out}(t) A_{ij} - \sum_{j \in \delta_p^-(i)}\varphi_{ji}^{in}(t) A_{ji}  \notag & \\ 
& ~\cdots+ \sum_{j\in \delta_c^+(i)} f_{ij} - \sum_{j \in \delta_c^-(i)} f_{ji} \label{eq:nodal-balance-j} & \\
& ~\cdots= \sum_{r \in \mathcal R_i} f_r(t) - \sum_{d \in \mathcal D_i} f_d(t) - \sum_{s \in \mathcal S_i} f_s(t). &
 \notag
\end{flalign}
For any storage $s \in \mathcal S$, and for every junction $i$ in the well corresponding to that storage, we have
\begin{subequations}
\begin{flalign}
& f_s(t) = f_s^{wh}(t), & \\ 
& \frac{f_s^{wh}(t)}{A_w} = \sum_{j \in \delta_p^-(i)} \varphi_{ji}^{in}(t) \text{ if $i$ is a well head}, \label{eq:wh-nodal-balance} & \\ 
& \sum_{j \in \delta_p^+(i)} \varphi_{ij}^{out}(t)  - \sum_{j \in \delta_p^-(i)} \varphi_{ji}^{in}(t)  = 0, \label{eq:intermediate-nodal-balance} & \\ 
& \frac{f_s^{bh}(t)}{A_w} = \sum_{j \in \delta_p^+(i)} \varphi_{ij}^{out}(t) A_{ij} \text{ if $i$ is a bottom hole,} \label{eq:bh-nodal-balance} &
\end{flalign}
\label{eq:nodal-balance-w}
\end{subequations}
where $A_w$ is the cross-sectional area of the well. 

\subsection{Objective Function} \label{subsec:obj} 
Thus far, we have described all the variables and equations governing the operation of a  pipeline network. Formally, we define day-ahead operation by specifying $\alpha_{ij},f_r(t), f_d(t), f_s(t)$ in Eq. \eqref{eq:compressor_var}, \eqref{eq:receipt_var}, \eqref{eq:delivery_var},  and \eqref{eq:storage_var} respectively. 
The objective function is intended to optimize the business goals of the pipeline operator, as well as the operating efficiency of the pipeline system system as a whole. The business goal of pipeline operators is to maximize the economic value of transportation service, given by
\begin{flalign} \label{eq:obj_profit}
J_P\ & \triangleq  \sum_{d \in \mathcal D} \int_0^T c_{d}(t)f_d(t)~ dt + \sum_{r \in \mathcal R} \int_0^T c_{r}(t)f_r(t)~ dt.
\end{flalign}
Here, $c_{d}(t)\geqslant 0$ is the bid (or offer price) of the consumer at off-take point $d \in \mathcal D$, and $c_{r}(t)\leqslant 0$ is the bid of the supplier at gas intake point $r \in \mathcal R$.

To promote efficiency, the system operators minimize the power used by compressors, which is given by
\begin{flalign} \label{eq:obj_eff}
J_E  & \triangleq  \sum_{(i,j)\in\mathcal{C}} \int_0^T W_{ij}(t) \cdot f_{ij}(t) \cdot \varphi_0 ~dt.
\end{flalign}
The overall objective function of the optimal control problem is a linear combination of $J_P$ and $J_E$, with a weight of $\kappa \in [0, 1]$ on the economic value, and is stated for minimization as
\begin{flalign} \label{eq:obj}
J \triangleq \kappa (-J_P) + (1-\kappa) J_E, 
\end{flalign}
The limits $\kappa = 0$ and $\kappa=1$  correspond to the cases when the objective is either efficiency or economic value. The entire optimal control problem can be formulated as
\begin{equation} \label{prob:ocp1}
\begin{array}{ll}
\!\!\!\! \mathrm{min}  & J \text{ in } \eqref{eq:obj} \\
\!\!\!\! \text{s.t.} & \text{Pipe dynamics } \eqref{eq:pipe-ode-horz}, \eqref{eq:in}, \text{ and } \eqref{eq:out}, \eqref{eq:pipe-ode-horz}, \\
& \text{Well dynamics } \eqref{eq:in}, \eqref{eq:out}, \eqref{eq:pipe-ode-1} \text{ with } g_{\parallel} = g, \\
& \text{Compressor action } \eqref{eq:c-1} \text{ or } \eqref{eq:c-2}, \eqref{eq:work-compressor}, \text{ \& } \eqref{eq:boost}, \\
& \text{Nodal \& well flow balance } \eqref{eq:nodal-balance-j}, \eqref{eq:wh-nodal-balance}-\eqref{eq:bh-nodal-balance}, \\
& \text{Pipeline pressure \& flow limits } \eqref{eq:rho_junction_limits}-\eqref{eq:pipe_flux_limits},  \\
& \text{Compressor limits } \eqref{eq:compressor_flow_limits}-\eqref{eq:compressor_power_limits}, \\
& \text{Storage facility limits } \eqref{eq:storage_limits}-\eqref{eq:cr_limits},\\
& \text{Receipt and delivery limits } \eqref{eq:receipt_limits}-\eqref{eq:delivery_limits}, \\
& \text{Time periodicity } \eqref{eq:endpoint_cond}.
\end{array}
\end{equation}

\section{Computational Results} \label{sec:results} 
We now present the results of several computational experiments to illustrate the function of the optimal control problem defined above. The problem \eqref{prob:ocp1} is solved via a primal-dual interior point method after approximating the derivatives in the component physics constraints using finite differences, and approximating integrals using the trapezoidal rule. We use a coarse (one hour) time discretization, which enables tractability of large-scale problems of interest.  The use of such coarse discretization has been validated in the regime of interest in several previous studies \cite{Zlotnik2015,Sundar2018,Sundar2019,Zlotnik2019,Roald2020}.  The mathematical formulation was implemented using the Julia programming language \cite{Julia-2017} with JuMP v0.21.3 \cite{Dunning2017} and KNITRO \cite{Byrd2006} as the primal-dual interior point solver. The computational experiments were run on a MacBook Pro with a 2.9 GHz dual-core processor and 16 GB RAM. The optimization time horizon is set to $T = 24$ hours and the value of $\kappa$ in Eq. \eqref{eq:obj} is set to $0.95$.

\subsection{Description of the Test Cases} \label{subsec:test-cases}
We consider two test cases: i) a $6$-junction network, and ii) a larger test network synthesized based on the Questar pipeline network located in the central United States (US). 
\subsubsection{6-junction network without storage} \label{subsubsec:no-storage-6}
\begin{figure}[!ht]
    \centering
    \includegraphics[scale=0.8]{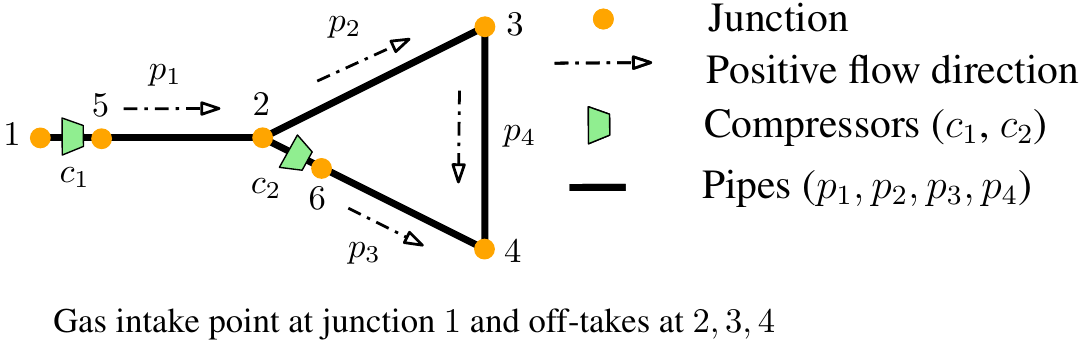}
    \caption{6-junction network without storage.}
    \label{fig:6-node-schematic}
\end{figure}

The $6$-junction network without storage contains a single loop, two compressors, \edit{$5$ gas off-take points with ids $d_1, d_2, \dots, d_5$ and one gas intake-point $s_1$. The off-take point ids and the junction at which they are located is given in Table \ref{tab:off-take-6-junction}.
\begin{table}[!ht]
    \centering
    \caption{The id of the junction where each off-take point is co-located for the $6$-junction network.}
    \begin{tabular}{|c|c|c|c|c|c|}
        \hline
        gas off-take point id & $d_1$ & $d_2$ & $d_3$ & $d_4$ & $d_5$ \\
        \hline
        junction id & 2 & 3 & 4 & 3 & 4 \\
        \hline
    \end{tabular}
    \label{tab:off-take-6-junction}
\end{table}
Throughout the rest of the article, we use $i$ and $d_i$ interchangeably to refer to the off-take point $d_i$. The one gas intake point (supply), located at junction $1$, which is a slack junction with pressure (or correspondingly, density) fixed at $4$ \si{\mega\pascal} or $580$ \si{psi}. The minimum and maximum operating pressures in the entire pipeline network are set to $3$ \si{\mega\pascal} and $6$ \si{\mega\pascal}, respectively. The supply and demand prices of gas at the intake and off-take points are constants for the entire time horizon of interest. The supply price of $1$ \si{\kilogram\per\second} of gas at the intake point is set to -$1.24$ dollars (the negative sign indicates payments to suppliers) and the sale price of $1$ \si{\kilogram\per\second} of gas at the off-take points $1$, $2$, $3$, $4$, and $5$ are set to $3$, $4$, $5$, $2.5$, and $3$ dollars, respectively. The maximum gas withdrawal rate at the off-take points is a time-varying function $f_d^{\max}(t)$, as shown in Fig. \ref{fig:baseline-results} (a). The maximum rate of gas injection into the network at an intake-point is designated by $f_r^{\max}(t)$, and is set to a constant value of $150$ \si{\kilogram\per\second} for all $t$. The friction factor for all the pipes was set to a value of $0.01$ and the length and diameter of the $4$ pipes are given in Table \ref{tab:ld}.}
\begin{table}[!ht]
    \centering
    \caption{Length and diameter of the pipes in the $6$-junction network.}
    \begin{tabular}{|c|c|c|c|c|}
        \hline
        pipe id & $p_1$ & $p_2$ & $p_3$ & $p_4$ \\
        \hline
        length (\si{\km}) & 50 & 80 & 80 & 80 \\
        \hline
        diameter (\si{\metre}) & 0.6 & 0.6 & 0.6 & 0.3\\
        \hline
    \end{tabular}
    \label{tab:ld}
\end{table}
\edit{All the other operating limits and time-varying parameter data associated with the network are available in \url{https://github.com/lanl-ansi/GasModels.jl/tree/master/test/data}. The $6$-junction network (without any storage units) is shown in Fig. \ref{fig:6-node-schematic}.}
\begin{figure*}[!ht]
    \centering
    \includegraphics[scale=1.1]{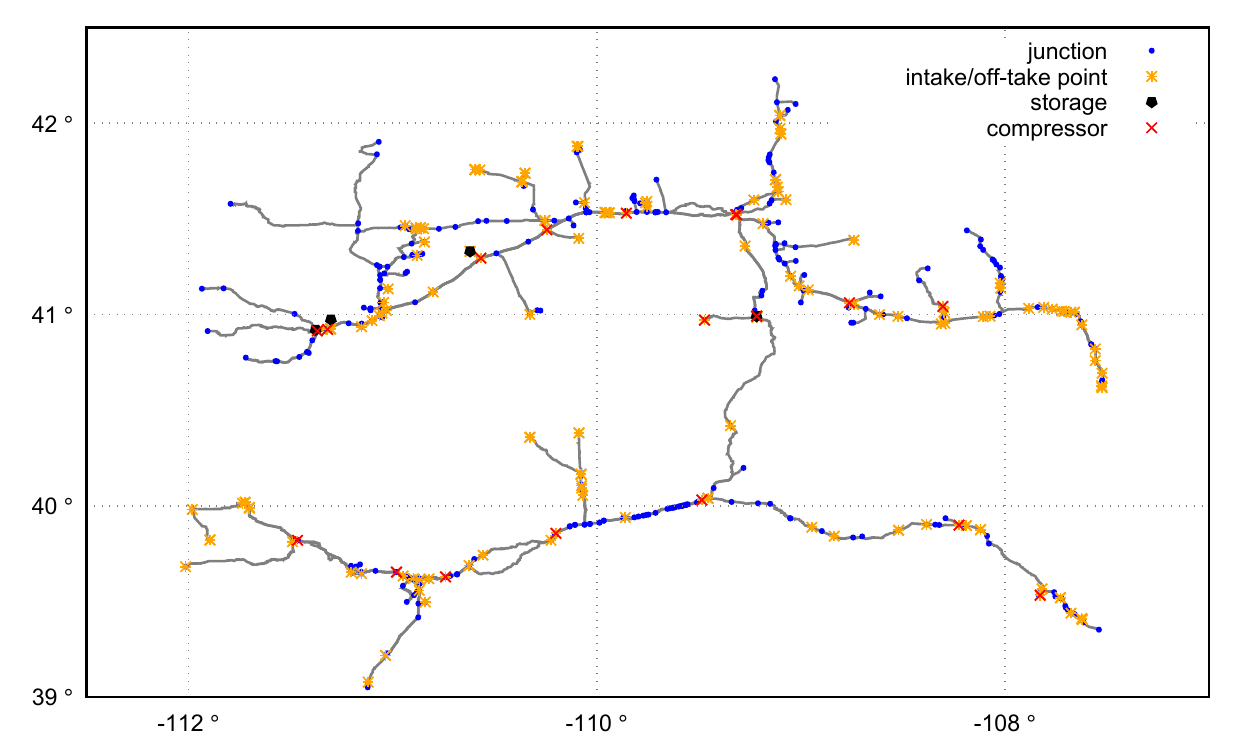}
    \caption{Schematic of the Questar pipeline network}
    \label{fig:questar}
\end{figure*}

\subsubsection{6-junction network with storage} \label{subsubsec:with-storage-6}
This network has a storage facility added to junction $3$ in the $6$-junction network. We assume that the when the reservoir is filled to its capacity, the pressure is $1365$ \si{psi}, and that the reservoir is $80$\% full at the beginning of the optimization time horizon. The minimum and maximum mass of gas in the reservoir is set to $3.5\times 10^8$ and $6.20\times 10^8$ \si{\kilogram} respectively. This reflects the typical design in which the reservoir is sustains continuous gas withdrawal from the storage for one to two months.

\subsubsection{Questar pipeline network} \label{subsubsec:questar}
The Questar network is an interstate natural gas pipeline in the US that provides transportation and underground storage services in the states of Utah, Wyoming, and Colorado. We synthesize a model from openly available data that has a total of $506$ junctions, $20$ compressors, $4$ natural gas storage facilities, $196$ intake/off-take points, and a total pipeline length of $3490$ \si{\km}. The schematic of the Questar network is shown in Fig. \ref{fig:questar}. The goal is to demonstrate the scalability of our approach, rather that to claim a particular analysis of the actual Questar system.\\

\subsection{Solution for the 6-junction network without storage} \label{subsec:baseline-6} 
\begin{figure}[!ht]
    \centering
    \includegraphics[scale=1.1]{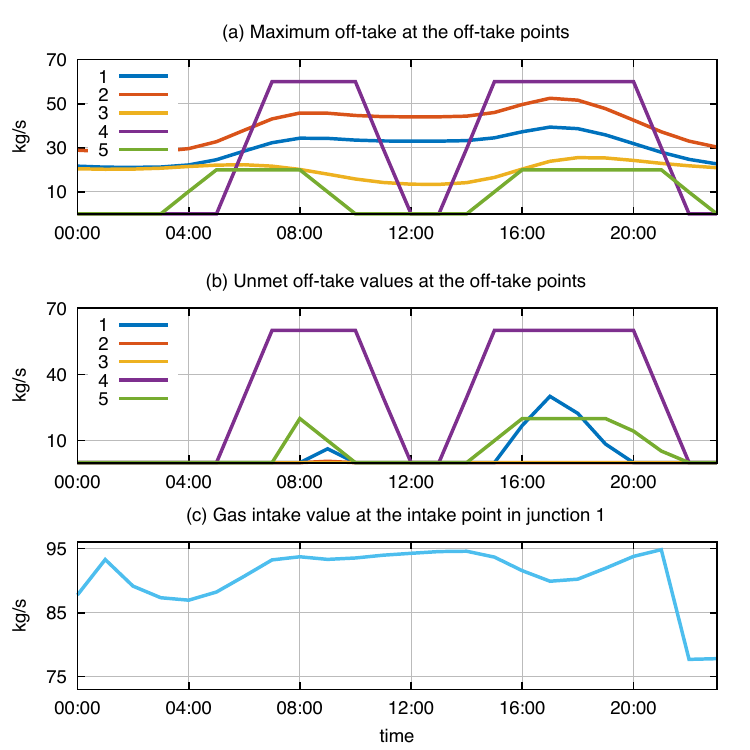}
    \caption{Results of the optimal control problem on the $6$-junction network without storage. Top: Maximum off-takes at the off-take points; Center: Off-take curtailment (unmet demand) at the off-take points; Bottom: Gas intake at the intake point at junction 1. The gas delivered at the off-take points is driven by the sale price of natural gas at that demand point. In this case, the off-take point $4$ does not receive any gas as the sale price at that point was the least. The sale prices at off-take points $3$  and $2$ were the two highest values, and hence, the objective function of the optimal control problem ensures that it delivers the maximum amount of gas at these two off-take points}
    \label{fig:baseline-results}
\end{figure}
Our first set of computational experiments aims to analyze the solution of the optimal control problem \eqref{prob:ocp1} for the $6$-junction pipeline network without any storage facilities.
We refer to the solution obtained by solving the optimal control problem for this case without storage as the baseline solution. The amount of gas that goes undelivered at the off-take points for the baseline solution is shown in Fig. \ref{fig:baseline-results} (b). The gas delivered at the off-take points is driven by the sale price of natural gas at that demand point. In this case, the off-take point $4$ does not receive any gas as the sale price at that point is smallest. The sale prices at off-take points $3$ and $2$ were the highest and second highest, so the optimal solution will attempt to deliver the maximum amount of gas at these two off-take points (also shown in the Fig. \ref{fig:baseline-results} (b)). Finally, Fig. \ref{fig:baseline-results} (c) shows the amount of gas that is entering the system from the intake point in junction $1$ of the network.

\subsection{Solution for the 6-junction network with storage} \label{subsec:bws-6}
Next we solve the optimal control problem \eqref{prob:ocp1} for the $6$-junction network with storage, and compare the solution with respect to the baseline solution of the previous section. 
\begin{figure}[!ht]
    \centering
    \includegraphics[scale=1.1]{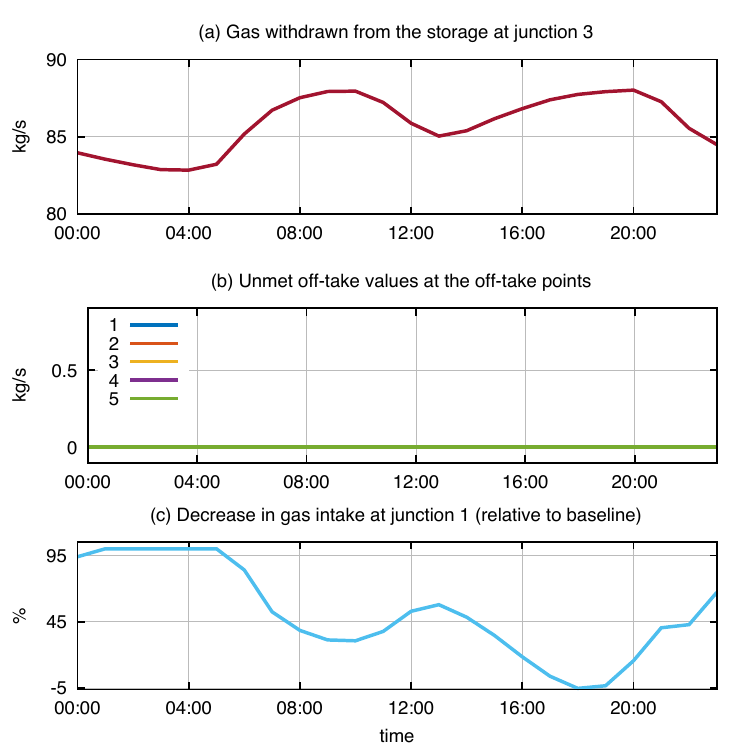}
    \caption{Results of the optimal control problem on the $6$-junction network with storage. Top: Maximum off-takes at the off-take points; Center: Off-take curtailment (unmet demand) at the off-take points; Bottom: Gas intake at the intake point at junction 1.  Notice that storage is able to support the baseline capacity and ensure that demands are met completely.}
    \label{fig:baseline-results-1}
\end{figure}
The results are presented in Fig. \ref{fig:baseline-results-1} (a), which shows the amount of gas that is withdrawn from storage for the full operating time horizon. The gas withdrawal from the storage causes a decrease in the amount of gas that is taken into the system from the intake point in junction $1$ (as high as a $96$\% decrease), as shown in Fig. \ref{fig:baseline-results-1} (c). The reason for this behaviour is that gas taken from the storage is not priced in the objective function in Eq. \eqref{eq:obj}, whereas gas from the intake point has a positive value for buy price. Furthermore, as observed in Fig. \ref{fig:baseline-results-1}, we see that unlike the baseline solution, all of the customers' requirements in the off-take points have been satisfied using the gas from storage. The next section presents results that illustrate situations where the storage can act both as a source and a sink within a $24$-hour time period. 

\begin{figure}[!ht]
    \centering
    \includegraphics[scale=1.1]{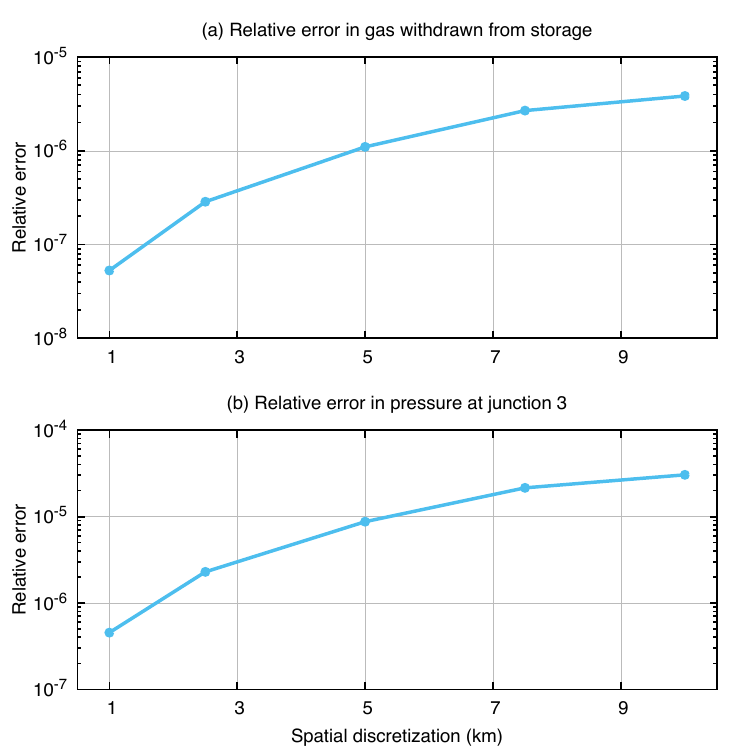}
    \caption{Plot of the $E_{s}(\delta)$ and $E_p(\delta)$ values for $\delta \in \{1, 2.5, 5.0, 7.5, 10.0\}$.}
    \label{fig:mesh-independence}
\end{figure}

\edit{\subsection{Spatial Discretization Study} \label{subsec:mesh} 
Here we provide a computational study to examine the sensitivity of the solution quality to the choice of spatial discretization used on the 6-junction network with storage.  The spatial discretization for the pipes in all the test cases is set to a value of $10$ \si{\km}. Thus, we vary the spatial discretization for the pipes in the network and examine the errors in the gas intake value at the storage facility and the pressure of the junction where the storage facility is located (junction $3$) relative to the smallest discretization value. The transient optimal control problem in Eq. \eqref{prob:ocp1} is solved with spatial discretization values in the set $\{0.5, 1, 2.5, 5.0, 7.5, 10.0\}$ \si{\km}. For each of these values, the metric we use to assess solution quality is the mean relative error between the gas intake from the storage and the pressure at junction $3$ and the same quantities computed using a very small spatial discretization of $0.5$ \si{\km} over the entire optimization time horizon. Mathematically, if we let $\delta$, $f_s(t, \delta)$, and $p_3(t, \delta)$ denote the spatial discretization, the gas intake from storage for the spatial discretization value of $\delta$, and the pressure at junction $3$ for the spatial discretization value of $\delta$, then, the relative errors in gas intake and pressure are defined by the following equations:}
\begin{subequations}
\begin{flalign}
& E_{s}(\delta) = \frac 1T \int_0^T  \frac{ \left|f_s(t, \delta) - f_s(t, 0.5)\right|}{ f_s(t, 0.5)} ~dt\label{eq:storage-error} \\
& E_{p}(\delta) = \frac 1T \int_0^T  \frac{ \left|p_3(t, \delta) - p_3(t, 0.5)\right|}{ p_3(t, 0.5)} ~dt \label{eq:pressure-error}
\end{flalign}
\label{eq:rel-errors}
\end{subequations}
\edit{The values of the relative errors for each $\delta \in \{1, 2.5, 5.0, 7.5, 10.0\}$ are shown using a logarithmic axis in Fig. \ref{fig:mesh-independence}. As observed from the Fig. \ref{fig:mesh-independence}, the mean relative error in the both the intake and pressure values is less than $10^{-4}$ even for a spatial discretization value of $10$ \si{\km}, which we believe justifies the use of this value for the rest of the computational experiments.}

\subsection{When storage alternates between injection and withdrawal} \label{subsec:with-storage-lb-6}
\begin{figure}[ht]
    \centering
    \includegraphics[scale=1.1]{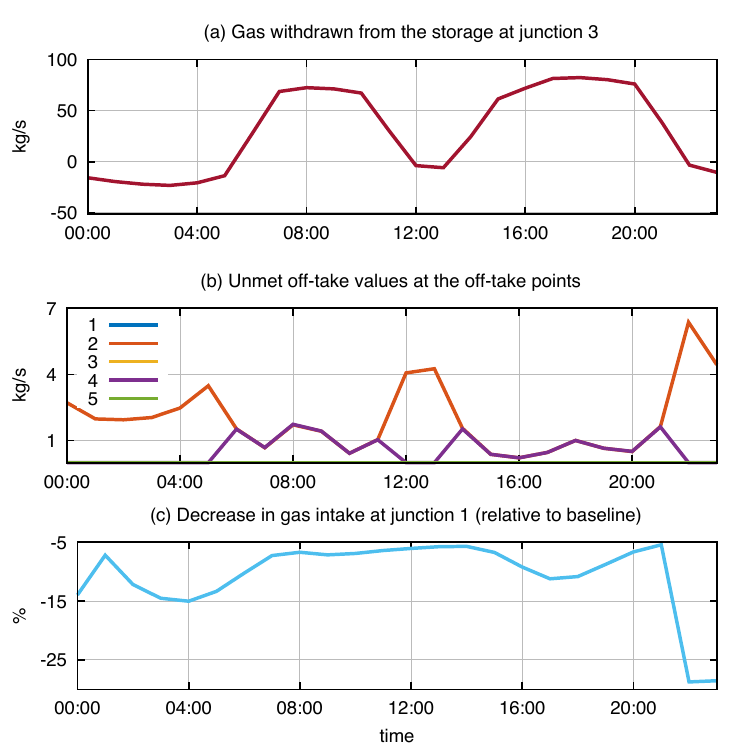}
    \caption{Results of the optimal control problem on the $6$-junction network with a lower bound on the intake point in junction $1$. Top: Maximum off-takes at the off-take points; Center: Off-take curtailment (unmet demand) at the off-take points; Bottom: Gas intake at the intake point at junction 1. If there is gas surplus to demand, it is injected into storage; else the storage permits withdrawal to meet the demand.}
    \label{fig:baseline-results-2}
\end{figure}
Here we show that if a storage facility has the required equipment, it can be beneficial to operate in both injection and withdrawal phases during a $24$-hour time horizon.  To do this, we increase the lower bound on the intake to $100$ \si{\kilogram\per\second}, so that in the optimal solution the network draws a minimum of $100$ \si{\kilogram\per\second} from the intake point so that the surplus gas is then  stored temporarily in the storage facility before it gets used to satisfy the demand of customers at the off-take points. The solution obtained for this case is shown in Fig. \ref{fig:baseline-results-2}. Here, the optimal solution dictates that gas is injected into the storage for the first four hours and again near the twelfth and the twenty-third hours of the operating horizon, and  is withdrawn during the rest of the optimization interval. In Fig. \ref{fig:baseline-results-2} (a), negative values indicate gas is being injected into the storage facility. Also, unlike the results in Sec. \ref{subsubsec:with-storage-6}, not all customers' requirement at the off-take points are satisfied fully as seen from  Fig. \ref{fig:baseline-results-2} (b). Similar to the results in the previous section, the decrease in the amount of gas taken in from the intake point, relative to the baseline solution in Sec. \ref{subsec:baseline-6}, is shown in Fig. \ref{fig:baseline-results-2} (c).

\subsection{Value of Storage Modeling} \label{subsec:storage-results}
\begin{figure}[!ht]
    \centering
    \includegraphics[scale=1.1]{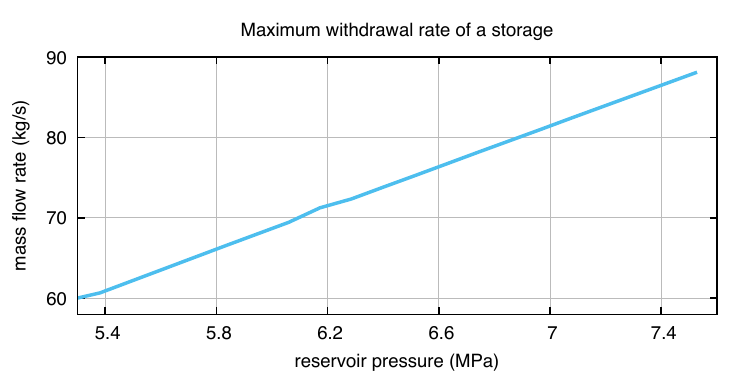}
    \caption{Reservoir Pressure vs maximum withdrawal rate of the storage facility. As the plot indicates, the maximum amount of gas that can be withdrawn from the storage facility decreases as more gas is being withdrawn and the pressure falls. The decrease reflects the physics of the reservoir as a volume with finite capacity, in sharp contrast to a battery-like model of storage, wherein, the maximum withdrawal capacity would have remained constant.}
    \label{fig:storage-rate}
\end{figure}

We previously discussed how the only storage facility models that have been previously integrated into pipeline operational optimization problems represent reservoirs as simple batteries with capacity limits \cite{He2017}. Such models do not include the actual physics  of  storage  facilities, but rather are akin to battery storage units in electric power networks where the maximum amount of power drawn from a battery does not vary based on state of charge. The maximum withdrawal directly reduces to zero from its rated power limit once the battery is completely discharged.
\emph{One of the main contributions of this work is a model of storage that is governed by the physics, and is yet amenable to integration into an optimal control framework for natural gas operations. We now illustrate qualitatively why battery-like behaviour in the storage model is not physically realistic.} 
For this study, we examine a single storage facility in isolation. We fix the pressure at the well-head to its minimum value of 250 \si{psi}, vary the reservoir pressure between its limits, and compute the amount of gas that can be withdrawn from the storage facility through the well for each value of the reservoir pressure. This gives the maximum rate at which gas can be withdrawn from the storage for a specified value of reservoir pressure. The plot of this maximum withdrawal rate as a function of the reservoir pressure is shown in Fig. \ref{fig:storage-rate}. As the plot indicates, the maximum amount of gas that can be withdrawn from the storage facility decreases as more gas is being withdrawn. The decrease reflects the physics of the reservoir as a volume with finite capacity, in sharp contrast to a battery-like model of storage.  

\subsection{Scalability study} \label{subsec:scalability}
\begin{figure}
    \centering
    \includegraphics[scale=1.1]{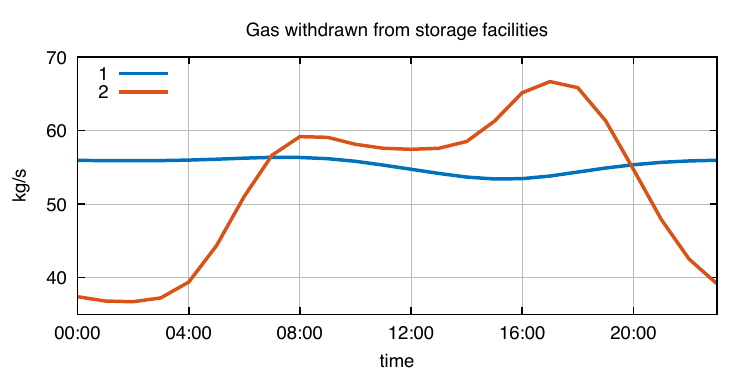}
    \caption{Mass flow rate of the gas injected into the Questar pipeline network from the two storage facilities. The results of the optimal control problem on the Questar pipeline network indicated that the other two storage facilities were neither injecting or withdrawing gas from the pipeline. }
    \label{fig:questar-storage}
\end{figure}

The final set of results presented corroborate the effectiveness and scalability  of the reduced order models and the resulting optimal control problem developed in this article by studying a model based on the Questar network, a real pipeline system in the US. The transient optimal control problem required 3 minutes of computation time to obtain the optimal solution. The Questar network has $4$ storage units; the solution to the optimal control problem had two storage units injecting gas into the pipeline network, and the remaining two storage were not active. The Fig. \ref{fig:questar-storage} shows the amount of gas that is being injected into the Questar pipeline network from the two storage facilities. 

\section{Conclusion} \label{sec:conclusion}

We have presented a new nonlinear optimal control problem for intra-day operation of a natural gas pipeline network, which incorporates unprecedented physics-based reduced order modeling of storage reservoirs for which wellhead flow rates depend on the inventory as well as bottom hole and wellhead pressures.  Expanding on studies that developed pipeline transient optimization methods, we maximize an objective function that quantifies economic profit and network efficiency subject to physical flow equations as well as engineering and operating limitations.  The problem is solved using a primal-dual interior point solver, and performance is validated in computational experiments and simulations on a small test network as well as a large and realistic test network based on a real pipeline.  Many pipeline systems include underground storage facilities in their original design in order to service seasonal demand.  Continued and efficient operation of these facilities is often required for these pipelines to function securely at their designed capacity.   The complex relation between pressure gradients and flows in pipeline systems has precluded storage facilities and flowing supplies from being utilized at their most efficient capacities simultaneously to service time-dependent loads.  The optimal control methodology for jointly scheduling pipeline flows with storage facility operations provides a promising path to address this critical challenge.

\section*{Acknowledgements}
This work was supported by the U.S. Department of Energy's Advanced Grid Modeling (AGM) projects \emph{Joint Power System and Natural Gas Pipeline Optimal Expansion} and \emph{Dynamical Modeling, Estimation, and Optimal Control of Electrical Grid-Natural Gas Transmission Systems}. The research work conducted at Los Alamos National Laboratory is done under the auspices of the National Nuclear Security Administration of the U.S. Department of Energy under Contract No. 89233218CNA000001. We gratefully thank the AGM program managers Alireza Ghassemian and Sandra Jenkins for their support of this research. 

\bibliographystyle{IEEEtran}
\bibliography{references.bib}

\begin{thebibliography}{10}
\providecommand{\url}[1]{#1}
\csname url@samestyle\endcsname
\providecommand{\newblock}{\relax}
\providecommand{\bibinfo}[2]{#2}
\providecommand{\BIBentrySTDinterwordspacing}{\spaceskip=0pt\relax}
\providecommand{\BIBentryALTinterwordstretchfactor}{4}
\providecommand{\BIBentryALTinterwordspacing}{\spaceskip=\fontdimen2\font plus
\BIBentryALTinterwordstretchfactor\fontdimen3\font minus
  \fontdimen4\font\relax}
\providecommand{\BIBforeignlanguage}[2]{{%
\expandafter\ifx\csname l@#1\endcsname\relax
\typeout{** WARNING: IEEEtran.bst: No hyphenation pattern has been}%
\typeout{** loaded for the language `#1'. Using the pattern for}%
\typeout{** the default language instead.}%
\else
\language=\csname l@#1\endcsname
\fi
#2}}
\providecommand{\BIBdecl}{\relax}
\BIBdecl

\bibitem{EIA-1}
``{U.S.} natural gas consumption sets new record in 2019,''
  \url{https://www.eia.gov/todayinenergy/detail.php?id=43035}.

\bibitem{EIA-2}
``Competition between coal and natural gas affects power markets,''
  \url{https://www.eia.gov/todayinenergy/detail.php?id=31672}.

\bibitem{Rinaldi2001}
S.~M. Rinaldi, J.~P. Peerenboom, and T.~K. Kelly, ``Identifying, understanding,
  and analyzing critical infrastructure interdependencies,'' \emph{IEEE control
  systems magazine}, vol.~21, no.~6, pp. 11--25, 2001.

\bibitem{Li2008}
T.~Li, M.~Eremia, and M.~Shahidehpour, ``Interdependency of natural gas network
  and power system security,'' \emph{IEEE Transactions on Power Systems},
  vol.~23, no.~4, pp. 1817--1824, 2008.

\bibitem{Roald2020}
L.~A. {Roald}, K.~{Sundar}, A.~{Zlotnik}, S.~{Misra}, and G.~{Andersson}, ``An
  uncertainty management framework for integrated gas-electric energy
  systems,'' \emph{Proceedings of the IEEE}, vol. 108, no.~9, pp. 1518--1540,
  2020.

\bibitem{Zlotnik2015}
A.~Zlotnik, M.~Chertkov, and S.~Backhaus, ``Optimal control of transient flow
  in natural gas networks,'' in \emph{2015 54th IEEE conference on decision and
  control (CDC)}.\hskip 1em plus 0.5em minus 0.4em\relax IEEE, 2015, pp.
  4563--4570.

\bibitem{Sundar2018}
K.~Sundar and A.~Zlotnik, ``State and parameter estimation for natural gas
  pipeline networks using transient state data,'' \emph{IEEE Transactions on
  Control Systems Technology}, vol.~27, no.~5, pp. 2110--2124, 2018.

\bibitem{Sundar2019}
------, ``Dynamic state and parameter estimation for natural gas networks using
  real pipeline,'' in \emph{2019 IEEE Conference on Control Technology and
  Applications (CCTA)}.\hskip 1em plus 0.5em minus 0.4em\relax IEEE, 2019, pp.
  106--111.

\bibitem{Jalving2018}
J.~Jalving and V.~M. Zavala, ``An optimization-based state estimation framework
  for large-scale natural gas networks,'' \emph{Industrial \& Engineering
  Chemistry Research}, vol.~57, no.~17, pp. 5966--5979, 2018.

\bibitem{Zlotnik2019}
A.~Zlotnik, K.~Sundar, A.~M. Rudkevich, A.~Beylin, and X.~Li, ``Optimal control
  for scheduling and pricing intra-day natural gas transport on pipeline
  networks,'' in \emph{2019 IEEE 58th Conference on Decision and Control
  (CDC)}.\hskip 1em plus 0.5em minus 0.4em\relax IEEE, 2019, pp. 4887--4884.

\bibitem{wylie1978fluid}
E.~B. Wylie and V.~L. Streeter, ``Fluid transients,'' \emph{mhi}, 1978.

\bibitem{Osiadacz1984}
A.~Osiadacz, ``Simulation of transient gas flows in networks,''
  \emph{International journal for numerical methods in fluids}, vol.~4, no.~1,
  pp. 13--24, 1984.

\bibitem{misra2020monotonicity}
S.~Misra, M.~Vuffray, and A.~Zlotnik, ``Monotonicity properties of physical
  network flows and application to robust optimal allocation,''
  \emph{Proceedings of the IEEE}, vol. 108, no.~9, pp. 1558--1579, 2020.

\bibitem{Chapman2005}
K.~S. Chapman, P.~Krishniswami, V.~Wallentine, M.~Abbaspour, R.~Ranganathan,
  R.~Addanki, J.~Sengupta, and L.~Chen, ``Virtual pipeline system testbed to
  optimize the us natural gas transmission pipeline system,'' Kansas State
  University, Tech. Rep., 2005.

\bibitem{Grundel2013}
S.~Grundel, N.~Hornung, B.~Klaassen, P.~Benner, and T.~Clees, ``Computing
  surrogates for gas network simulation using model order reduction,'' in
  \emph{Surrogate-Based Modeling and Optimization}.\hskip 1em plus 0.5em minus
  0.4em\relax Springer, 2013, pp. 189--212.

\bibitem{Zlotnik2015a}
A.~Zlotnik, S.~Dyachenko, S.~Backhaus, and M.~Chertkov, ``Model reduction and
  optimization of natural gas pipeline dynamics,'' in \emph{Dynamic Systems and
  Control Conference}, vol. 57267.\hskip 1em plus 0.5em minus 0.4em\relax
  American Society of Mechanical Engineers, 2015, p. V003T39A002.

\bibitem{Zlotnik2016}
A.~Zlotnik, L.~Roald, S.~Backhaus, M.~Chertkov, and G.~Andersson, ``Coordinated
  scheduling for interdependent electric power and natural gas
  infrastructures,'' \emph{IEEE Transactions on Power Systems}, vol.~32, no.~1,
  pp. 600--610, 2016.

\bibitem{Reddy2006}
H.~P. Reddy, S.~Narasimhan, and S.~M. Bhallamudi, ``Simulation and state
  estimation of transient flow in gas pipeline networks using a transfer
  function model,'' \emph{Industrial \& engineering chemistry research},
  vol.~45, no.~11, pp. 3853--3863, 2006.

\bibitem{Alamian2012}
R.~Alamian, M.~Behbahani-Nejad, and A.~Ghanbarzadeh, ``A state space model for
  transient flow simulation in natural gas pipelines,'' \emph{Journal of
  Natural Gas Science and Engineering}, vol.~9, pp. 51--59, 2012.

\bibitem{He2017}
Y.~He, M.~Shahidehpour, Z.~Li, C.~Guo, and B.~Zhu, ``Robust constrained
  operation of integrated electricity-natural gas system considering
  distributed natural gas storage,'' \emph{IEEE Transactions on Sustainable
  Energy}, vol.~9, no.~3, pp. 1061--1071, 2017.

\bibitem{Basumatary2005}
R.~Basumatary, P.~Dutta, M.~Prasad, and K.~Srinivasan, ``Thermal modeling of
  activated carbon based adsorptive natural gas storage system,''
  \emph{Carbon}, vol.~43, no.~3, pp. 541--549, 2005.

\bibitem{Oldenburg2001}
C.~Oldenburg, K.~Pruess, and S.~M. Benson, ``Process modeling of co2 injection
  into natural gas reservoirs for carbon sequestration and enhanced gas
  recovery,'' \emph{Energy \& Fuels}, vol.~15, no.~2, pp. 293--298, 2001.

\bibitem{Lokhorst2005}
A.~Lokhorst and T.~Wildenborg, ``Introduction on co2 geological
  storage-classification of storage options,'' \emph{Oil \& gas science and
  technology}, vol.~60, no.~3, pp. 513--515, 2005.

\bibitem{Thorley1987}
A.~Thorley and C.~Tiley, ``Unsteady and transient flow of compressible fluids
  in pipelines—a review of theoretical and some experimental studies,''
  \emph{International journal of heat and fluid flow}, vol.~8, no.~1, pp.
  3--15, 1987.

\bibitem{Menon2005}
E.~S. Menon, \emph{Gas pipeline hydraulics}.\hskip 1em plus 0.5em minus
  0.4em\relax {CRC} Press, 2005.

\bibitem{Gyrya2019}
V.~Gyrya and A.~Zlotnik, ``An explicit staggered-grid method for numerical
  simulation of large-scale natural gas pipeline networks,'' \emph{Applied
  Mathematical Modelling}, vol.~65, pp. 34--51, 2019.

\bibitem{Dyachenko2017}
S.~A. Dyachenko, A.~Zlotnik, A.~O. Korotkevich, and M.~Chertkov, ``Operator
  splitting method for simulation of dynamic flows in natural gas pipeline
  networks,'' \emph{Physica D: Nonlinear Phenomena}, vol. 361, pp. 1--11, 2017.

\bibitem{Mak2016}
T.~W. Mak, P.~Van~Hentenryck, A.~Zlotnik, H.~Hijazi, and R.~Bent, ``Efficient
  dynamic compressor optimization in natural gas transmission systems,'' in
  \emph{2016 American Control Conference (ACC)}.\hskip 1em plus 0.5em minus
  0.4em\relax IEEE, 2016, pp. 7484--7491.

\bibitem{flanigan1995underground}
O.~Flanigan, \emph{Underground gas storage facilities: Design and
  implementation}.\hskip 1em plus 0.5em minus 0.4em\relax Elsevier, 1995.

\bibitem{jansen2008model}
J.-D. Jansen, O.~H. Bosgra, and P.~M. Van~den Hof, ``Model-based control of
  multiphase flow in subsurface oil reservoirs,'' \emph{Journal of Process
  Control}, vol.~18, no.~9, pp. 846--855, 2008.

\bibitem{chen2020frankenstein}
B.~Chen, D.~R. Harp, R.~J. Pawar, P.~H. Stauffer, H.~S. Viswanathan, and R.~S.
  Middleton, ``Frankenstein’s romster: Avoiding pitfalls of reduced-order
  model development,'' \emph{International Journal of Greenhouse Gas Control},
  vol.~93, p. 102892, 2020.

\bibitem{Julia-2017}
J.~Bezanson, A.~Edelman, S.~Karpinski, and V.~B. Shah, ``Julia: A fresh
  approach to numerical computing,'' \emph{SIAM {R}eview}, vol.~59, no.~1, pp.
  65--98, 2017.

\bibitem{Dunning2017}
I.~Dunning, J.~Huchette, and M.~Lubin, ``Jump: A modeling language for
  mathematical optimization,'' \emph{SIAM Review}, vol.~59, no.~2, pp.
  295--320, 2017.

\bibitem{Byrd2006}
R.~H. Byrd, J.~Nocedal, and R.~A. Waltz, ``{KNITRO}: An integrated package for
  nonlinear optimization,'' in \emph{Large-scale nonlinear optimization}.\hskip
  1em plus 0.5em minus 0.4em\relax Springer, 2006, pp. 35--59.

\end{thebibliography}

\end{document}